\documentclass[11pt]{article}
\usepackage{amsmath}
\usepackage{amssymb}
\usepackage{amsfonts}
\usepackage{mathrsfs}
\usepackage{shortvrb,psfrag}
\usepackage{enumerate}
\usepackage{color}
\usepackage[dvips]{graphicx}

\let\f=\frac

\let\Om=\Omega

\def\R{\Bbb R}

\def\no{\noindent}

\def\endproof{\hphantom{MM}\hfill\llap{$\square$}\goodbreak}

\newcommand{\beq}{\begin{equation}}
\newcommand{\eeq}{\end{equation}}
\newcommand{\ben}{\begin{eqnarray}}
\newcommand{\een}{\end{eqnarray}}
\newcommand{\beno}{\begin{eqnarray*}}
\newcommand{\eeno}{\end{eqnarray*}}

\newtheorem{Theorem}{Theorem}[section]
\newtheorem{Definition}[Theorem]{Definition}
\newtheorem{Proposition}[Theorem]{Proposition}
\newtheorem{Lemma}[Theorem]{Lemma}

\newtheorem{Remark}[Theorem]{Remark}

\begin{document}
\title{\bf Boundary Regularity Criteria for the 6D Steady Navier-Stokes and MHD Equations}

\author{Jitao LIU$^{\dag,\sharp}$, Wendong WANG$^{\ddag}$\\[2mm]
{\small $^{\dag}$ The Institute of Mathematical Sciences, The Chinese University of Hong Kong, Shatin, N.T., Hong Kong.}\\[1mm]
{\small $^{\ddag}$ School of  Mathematical Sciences, Dalian University of Technology, Dalian, 116024, P. R. China.}\\[1mm]
{\small $^{\sharp}$ Universidade Federal do Rio de Janeiro, Av. Athos da Silveira Ramos, 149, Ilha do Fund\~{a}o,}\\[1mm]
{\small  Rio de Janeiro, RJ 21941-909, Brazil.}\\[1mm]
{\small E-mail: lijt1005@gmail.com, wendong@dlut.edu.cn}}

\maketitle

\begin{abstract}
It is shown in this paper that suitable weak solutions to the 6D steady incompressible Navier-Stokes and MHD equations are H\"{o}lder continuous near boundary provided that either $r^{-3}\int_{B_r^+}|u(x)|^3dx$ or $r^{-2}\int_{B_r^+}|\nabla u(x)|^2dx$ is sufficiently small, which implies that the 2D Hausdorff measure of the set of singular points near the boundary is zero. This generalizes recent interior regularity results by Dong-Strain \cite{DS}.
\end{abstract}



\setcounter{equation}{0}
\section{Introduction}

In this paper, we consider the following 6D steady incompressible Navier-Stokes equations on $\Omega \subset \R^6$:
\begin{align} \label{eq:SNS}(\rm SNS)\,\, \left\{
\begin{aligned}
&-\Delta u+u\cdot \nabla u=-\nabla{\pi}+f,\\
&\quad\nabla\cdot u=0,\\
\end{aligned}
\right. \end{align}
where $u$ represents the fluid velocity field, $\pi$ is a scalar pressure, and  the boundary condition of $u$ is given as a no-slip condition, namely
\ben\label{bd}
u=0,&\hbox{on}&\partial\Omega.
\een
The analysis of the above equations is motivated by Struwe's question in \cite{St1}, where  the 5D steady Navier-Stokes equations was considered and he asked if analogous partial regularity results hold in spacial dimension $N>5.$ Recently interior regularity results in 6D are obtained by Dong-Strain \cite{DS}, and the main interests in this paper are in the boundary partial regularity for suitable weak solutions to the equations (\ref{eq:SNS}). As it's commented in \cite{DS}, six
is the highest dimension for stationary Navier-Stokes equations in which all the existing methods on partial regularity could be applied.

Recall the development of interior and boundary regularity criteria for the Navier-Stokes equations in brief. For the three dimensional time-dependent Navier-Stokes equations, partial regularity of weak solutions satisfying the local energy inequality was proved in a series of papers by Scheffer \cite{SV1,SV2,SV4}. Later, the notion of suitable weak solutions was first introduced in a celebrated paper by Caffarelli-Kohn-Nirenberg \cite{CKN}, showing that the set $\mathcal{S}$ of possible interior singular points of a suitable weak solution is one-dimensional parabolic Hausdorff measure zero. Simplified proofs and improvements can be found in many works by Lin \cite{Lin}, Ladyzhenskaya-Seregin \cite{LS}, Tian-Xin \cite{TX}, Escauriaza-Seregin-\v{S}ver\'{a}k \cite{ESS}, Seregin \cite{Se}, Gustafson-Kang-Tsai \cite{GKT}, Vasseur \cite{Va}, Kukavica \cite{Ku}, Wang-Zhang \cite{WZ2} and the references therein. Some similar boundary regularity results are also proved, see Seregin \cite{Se0, Se1, Se2}, Kang \cite{Ka, Ka0}, Gustafson-Kang-Tsai \cite{GKT0} and so on. Moreover, it's worth to mention that general curved boundary regularity was obtained by Seregin-Shilkin-Solonnikov in \cite{SSS}.

There are only fewer results available in the literature for the 4D and
higher dimensional time-dependent Navier-Stokes equations. In \cite{SV3}, Scheffer showed that
there exists a weak solution in $R^4\times R^+$, whose singular set has vanishing 3D Hausdorff measure.
Later, Dong-Du \cite{DD} proved that, for any
local-in-time smooth solution to the 4D Navier-Stokes equations, the
2D Hausdorff measure of the set of singular points at the first potential
blow-up time is equal to zero, and we refer to \cite{DG} for recent results  with general suitable weak solutions.

Now we turn to the steady Navier-Stokes equations. In a series of papers by Frehse and Ruzicka \cite{FR1,FR2,FR3,FR4}, the existence on a class of special regular solutions of (\ref{eq:SNS}) was obtained for the five-dimensional and higher dimensional case.  Gerhardt \cite{Ge} obtained the regularity of weak solutions under the  four-dimensional case. For $N\geq 5$, it is not known yet whether weak solutions are regular.
In \cite{St1,St2}, Struwe obtained partial regularity for $N=5$ by regularity methods of elliptic systems (c.f. Morrey \cite{Mo} and Giaqinta \cite{Gi}).
Later, the result of Struwe was extended to the boundary case by Kang \cite{Ka}. In a recent paper, Dong-Strain \cite{DS} extended the interior regularity result of Struwe to the six-dimensional space. Their main idea is to apply the iteration method and bootstrap arguments to get a suitable decay estimate of $L^{3/2}$ norm of $\nabla u$, then the Morrey lemma implies the required regularity.

In this paper, we generalized the result in \cite{DS} and proved the boundary regularity of six-dimensional steady Navier-Stokes equations. Moreover, we considered six-dimensional steady MHD equations, and obtained boundary regularity criteria independent of the magnetic field.  It is important to mention that the paper by H. Dong and X. Gu \cite{DG2} is available on the arXiv at almost the same time with ours, where they considered partial regularity of the 4D Navier-Stokes equations and 6D steady Navier-Stokes equations,  and Theorem  1.4 in \cite{DG2} is similar to the following Theorem \ref{thm:c} by a different approach. Compared with the interior case, the main difficulty lies in the low regularity of the pressure near boundary, the iterative scheme in \cite{DS} becomes invalid. To overcome the obstacle, the main idea in  \cite{DG2} is to firstly establish a weak decay estimate of certain scale-invariant quantities, and then successively improve
this decay estimate by a revised bootstrap argument; however, in our paper, we derive a new revised local energy inequality based on  a free parameter which plays a key role in the proof, and it allows us to obtain decay estimates of certain scale-invariant quantities.

At last, we refer to  \cite{FS} by Farwig and Sohr for existence and regularity criteria for weak solutions to inhomogeneous  Navier-Stokes equations.

Let us introduce the definition of suitable weak solutions near the boundary.
\begin{Definition} Let  $\Om\subset\mathbb{R}^6$ be an open domain, and $\Gamma\subset \partial\Omega$ be an open set. $(u,\pi)$ is said to be a suitable weak solution to the steady Navier-Stoks equations (\ref{eq:SNS}) in $\Om$ near the boundary $\Gamma$, if the following conditions hold.

(i)\,$u\in H^1(\Om),\,\pi\in L^{\f 32}(\Om),\,\nabla \pi\in L^{\f 65}(\Om),\,\, f\in L^{6}(\Om)$;

(ii)\,$(u,\pi)$ satisfies the equations(\ref{eq:SNS}) in the sense of distribution sense and the boundary condition $u|_{\Gamma}=0$ holds;

(iii)\,$u$ and $\pi$ satisfy the local energy inequality
\ben\label{eq:local energy}
2\int_{\Om}|\nabla u|^2\phi dx
\leq\int_{\Om}\big[|u|^2\triangle\phi+u\cdot\nabla\phi(|u|^2+2\pi)\big]+2 fu \phi dx
\een
for any nonnegative $C^{\infty}$ test function $\phi$ vanishing at the boundary $\partial\Omega\backslash\Gamma$ .
\end{Definition}

The existence of such a suitable weak solution can be found in \cite{FR3}. The major concern of this paper is the regularity and the main results can be stated as follows:

\begin{Theorem}\label{thm:c}
Let $(u,\pi)$ be a suitable weak solution to (\ref{eq:SNS}) in $B_1^+$ near the boundary $\{x\in B_1, x_6=0\}$.
Then $0$ is a regular point of $u$,
if there exists a small positive constant $\varepsilon$
such that one of the following conditions holds,
\beno
&&i) \limsup_{r\rightarrow0_{+}} r^{-3}\int_{B_r^+}|u(x)|^3dx<\varepsilon,\\
&&ii) \limsup_{r\rightarrow0_{+}}r^{-2}\int_{B_r^+}|\nabla u(x)|^2dx <\varepsilon.
\eeno
\end{Theorem}

\begin{Remark}
The  boundary regularity criteria above for the 6D steady Navier-Stokes equations generalize recent interior regularity results by Dong-Strain \cite{DS} and boundary regularity results for the 5D case by Kang in \cite{Ka0}. Here we consider the flat boundary for simplicity, the results hold true for general $C^2$ boundary, which follows from the analysis here and the standard techniques as in \cite{SSS}.
\end{Remark}

\begin{Theorem}\label{thm:h}
Let $(u,\pi)$ be a suitable weak solution to (\ref{eq:SNS}) in $B_1^+$ near the boundary  $\{x\in B_1, x_6=0\}$.
Then the $2D$ Hausdorff measure of the set of singular points of $(u,\pi)$ in $B_1^+$ is equal to zero.
\end{Theorem}

\begin{Remark}
This theorem follows directly from Theorem \ref{thm:c} by the standard arguments from the geometric measure theory, which is explained for example in \cite{CKN}.
\end{Remark}

As it will be shown later, Theorem \ref{thm:c} will follow from the following partial regularity criteria.

\begin{Proposition}\label{prop:cd}
Let $(u,\pi)$ be a suitable weak solution to (\ref{eq:SNS}) in $B_1^+$ near the boundary $\{x\in B_1, x_6=0\}$.
If there exists ${\rho}_0>0$ and a small positive constant $\varepsilon_1$
such that
\beno
\rho_0^{-3}\|u\|_{L^3(B_{\rho_0}^+)}^3+\rho_0^{-2}\|\nabla \pi\|_{L^{6/5}(B_{\rho_0}^+)}+\rho_0^3\|f\|_{L^3(B_{\rho_0}^+)}^3<\varepsilon_1
\eeno
Then $0$ is a regular point of $u$.
\end{Proposition}

\begin{Remark}
There are several remarks in order concerning the proof of the main results. First, it should be pointed out that it seems difficult to adapt the boundary regularity theory for 5D steady Navier-Stokes system by Kang in \cite{Ka0} to our case, since the key blow-up arguments there (Lemma 4.6, \cite{Ka0}) are based on the compact imbedding $W^{1,2}(B_1)\hookrightarrow L^3(B_1)$ which fails in the 6D case. Our analysis is motivated by a bootstrap argument due to Dong-Strain \cite{DS} for the interior regularity theory. However,  due to the boundary, new difficulties arise. In particular, there are slow decaying terms involving $E^{1/2}(\rho)$ in the pressure decomposition in the presence of boundaries (see Lemma \ref{lem:p1}),
where
$$E(\rho)=\rho^{-2}\int_{B_{\rho}^+}|\nabla u(x)|^2dx.$$
This means that $E^{1/2}(\rho)$ and $D_1(\rho)\equiv\rho^{-2}\parallel\nabla \pi\parallel_{L^{\f65}(B_{\rho}^+)}$ are the same order in the standard iterative scheme, which seems impossible to obtain an effective iterative estimate by the local energy inequality as in \cite{DS} (more details see Remark \ref{rmk:idea}).
To overcome the difficulty, we first derive a revised local energy inequality, see Proposition \ref{prop:local}, such that there exists a free parameter in the coefficient of the scaled energy on the large ball in the local energy inequality, which plays an important role in the required estimate in Lemma \ref{eq:iterative} and yields an effective iteration scheme.
\end{Remark}

As an application of the above results, we consider the 6D steady Magneto-hydrodynamics equations (MHD) as follows:
\begin{align} \label{eq:SMHD}(\rm SMHD)\,\, \left\{
\begin{aligned}
-\Delta u+u\cdot \nabla u&=-\nabla\pi+b\cdot \nabla b+f,\\
-\Delta b+u\cdot \nabla b&=b\cdot \nabla u ,\\
\nabla\cdot u&=0,\\
\nabla\cdot b&=0,
\end{aligned}
\right. \end{align}
where $u$, $b$ describe the fluid velocity field and the magnetic field
respectively, and $\pi$ is a scalar pressure. Similar to Navier-Stokes equations, the no-slip boundary condition is given as
\ben\label{eq:mbd}
u=b=0,&\hbox{on}&\partial\Omega.
\een
Now let us introduce the definition of suitable weak solutions of (\ref{eq:SMHD})-(\ref{eq:mbd}).
\begin{Definition} Let an open domain $\Om\subset\mathbb{R}^6$, and an open set $\Gamma\subset \partial\Omega$. We say that $(u,b,\pi)$ is a suitable weak solution of the steady MHD equations (\ref{eq:SMHD}) in $\Om$ near the boundary $\Gamma$, if the following conditions hold.

(i)\,\,$(u,b)\in H^1(\Om),\,\pi\in L^{\f 32}(\Om),\,\nabla \pi\in L^{\f 65}(\Om),\,\, f\in L^{6}(\Om)$;

(ii)\,\,$(u,b,\pi)$ satisfies the equations(\ref{eq:SMHD}) in distribution sense and the boundary condition $u|_{\Gamma}=b|_{\Gamma}=0$ holds;

(iii)\,\,$u$, $b$ and $\pi$ satisfy the local energy inequality
\ben\label{ieq:local energy inequlity-MHD}
&&2\int_{\Om}(|\nabla u|^2+|\nabla b|^2)\phi dx\nonumber\\
&\leq& \int_{\Om}\big[(|u|^2+|b|^2)\triangle\phi+u\cdot\nabla\phi(|u|^2+|b|^2+2\pi)-(b\cdot u)(b\cdot \nabla\phi)\big]+f\cdot u \phi dx,
\een
for any nonnegative $C^{\infty}$ test function $\phi$ vanishing at the boundary $\partial\Omega\backslash\Gamma$ .
\end{Definition}

Following the same route of Theorem \ref{thm:c} and using the idea of \cite{WZ1}, we can obtain boundary regularity criteria independent of the magnetic field.

\begin{Theorem}\label{thm:MHD}
Let $(u,b,\pi)$ be a suitable weak solution of (\ref{eq:SMHD}) in $B_1^+$ near the boundary $\{x\in B_1, x_6=0\}$.
Then $0$ is a regular point of $(u,b)$,
if there exists a small positive constant $\varepsilon$
such that one of the following conditions holds,
\beno
&&i)\,\,\limsup_{r\rightarrow0_{+}} r^{-3}\int_{B_r^+}|u(x)|^3dx<\varepsilon,\\
&&ii)\,\,\limsup_{r\rightarrow0_{+}}r^{-2}\int_{B_r^+}|\nabla u(x)|^2dx <\varepsilon.
\eeno
\end{Theorem}
\begin{Remark}
The above theorem can be seen as a generalization of Theorem \ref{thm:c} by assuming the magnetic field $b\equiv0.$ The novelty is that the boundary regularity criteria are independent of the magnetic field, and the reason comes from the absence of the term $\|b\|_3$ in the local energy inequality (\ref{ieq:local energy inequlity-MHD}).
\end{Remark}

The rest of the paper is organized as follows. In section 2, we introduce some notations and some technical lemmas about Stokes operator and local energy estimates; especially  a revised local energy inequality is obtained.
In section 3, we prove Theorem \ref{thm:c} under the assumption of Proposition \ref{prop:cd}.
Section 4 is devoted to the proof of Proposition \ref{prop:cd}. At last, we give a brief proof of Theorem \ref{thm:MHD}.

\section{Notations and some technical lemmas}

\setcounter{equation}{0}

Throughout this article, $C_0$ denotes an absolute constant independent of $u,\rho,r$
and may be different from line to line.

Let $(u,\pi)$ be a solution to the steady Navier-Stokes equations (\ref{eq:SNS}). Set the following scaling:
\ben\label{eq:scaling}
u^{\lambda}(x)={\lambda}u(\lambda x),\quad \pi^{\lambda}(x)={\lambda}^2\pi(\lambda x),\quad f^{\lambda}(x)={\lambda}^3f(\lambda x),
\een
for any $\lambda> 0,$ then the family $(u^{\lambda},\pi^{\lambda})$ is also a solution of (\ref{eq:SNS}) with $f$ replaced by $f^{\lambda}(x)$.
Now define some quantities which are invariant under the scaling (\ref{eq:scaling}):

$$A(r)=r^{-4}\int_{B_r^+}|u(x)|^2dx,\quad C(r)=r^{-3}\int_{B_r^+}|u(x)|^3dx;$$

$$E(r)=r^{-2}\int_{B_r^+}|\nabla u(x)|^2dx,\quad D_1(r)=r^{-2}\parallel\nabla \pi\parallel_{L^{\f65}(B_r^+)};$$

$$D(r)=r^{-3}\int_{B_r^+}|\pi-\pi_{B_r^+}|^{\f32}dx,\quad \pi_{B_r^+}=\frac{1}{|B_r^+|}\int_{B_r^+}\pi dx;$$

$$F(r)=r^{3}\int_{B_r^+}|f(x)|^3dx,$$
where
$B_r^+(x_0)$ is the  semi-ball of radius $r$ centered at $x_0$, and
we denote $B_r^+(0)$ by $B_r^+$.
Moreover, a solution $u$ is said to be regular at $x_0$ if $u\in L^\infty(B_r^+(x_0))$ for some $r>0$.

The following two lemmas about $L^p$ estimates on solutions to boundary value problems for the Stokes equations are well-known, for example, seeing
Solonnikov \cite{So}, Giga-Sohr \cite{GS}, Galdi \cite{Ga}, Seregin \cite{Se0} and Kang \cite{Ka}.

\begin{Lemma}\label{lem:stokes1}
Assume that $\Omega\subset \R^n$ is a bounded open domain , $n\geq 3$, $1<p<\infty$, and $f\in L^p(\Omega)$. $(v,q)$ satisfies the following conditions:
\begin{align} \label{eq:S}\,\, \left\{
\begin{aligned}
&v\in W^{2,p}(\Omega), \quad  q\in W^{1,p}(\Omega), \\
&-\triangle v+\nabla q=f,\quad \nabla\cdot v=0,\\
&\int_{\Omega}q dx=0,\quad v|_{\partial\Omega}=0.
\end{aligned}
\right. \end{align}
Then there holds
\beno
\|v\|_{W^{2,p}(\Omega)}+\|q\|_{W^{1,p}(\Omega)}\leq C(p,\Omega) \|f\|_{L^p(\Omega)}.
\eeno
\end{Lemma}

\begin{Lemma}\label{lem:stokes2}
Assume that $B_1^+\subset \R^n$ is a bounded open domain , $n\geq 3$, $1<p_0\leq p<\infty$, and $f\in L^p(B_1^+)$. $(v,q)$ satisfies the following conditions:
\begin{align} \label{eq:S}\,\, \left\{
\begin{aligned}
&v\in W^{2,p_0}(B_1^+), \quad q\in W^{1,p_0}(B_1^+), \\
&-\triangle v+\nabla q=f,\quad \nabla\cdot v=0,\\
&v|_{\{x_n=0\}}=0.
\end{aligned}
\right. \end{align}
Then there holds
\beno
\|v\|_{W^{2,p}(B_{1/2}^+)}+\|q\|_{W^{1,p}(B_{1/2}^+)}\leq C(p_0,p) \big(\|f\|_{L^p(B_1^+)}+\|q\|_{L^{p_0}(B_1^+)}+\|v\|_{W^{1,p_0}(B_1^+)}\big).
\eeno
\end{Lemma}

Recall the interior regular result by Dong-Strain (Proposition 3.7 and (3.33), \cite{DS}), which is necessary in the proof of Theorem \ref{thm:c}. We write it in a slightly different form.
\begin{Proposition}\label{prop:dong}
There exists $\varepsilon_0>0$ satisfying the following property.
Suppose that $(u,\pi)$ is a suitable weak solution of (\ref{eq:SNS}) in $B_1(x_0)$, and for some $\rho_0\in (0,1)$, it holds that
\beno
\rho_0^{-3}\int_{B_{\rho_0}(x_0)}|u|^3dx+ \rho_0^{-3}\int_{B_{\rho_0}(x_0)}|\pi-\pi_{B_{\rho_0}(x_0)}|^{\f32}dx+\rho_0^{3}\int_{B_{\rho_0}(x_0)}|f|^3dx\leq \varepsilon_2.
\eeno
Then, for $0<\rho<\rho_0/8$, the following inequality will hold uniformly
\ben\label{eq:dong}
\int_{B_{\rho}(x_0)}|\nabla u|^{\f65}dx\leq C_0 \rho^{\f{24}5+\f{2}{25}},
\een
where $C_0$ is a positive constant independent of $\rho$.
\end{Proposition}

Next, modifying the analysis in \cite{DS} by choosing a new test function, we will derive the following revised local energy
inequality, which improves the local energy
inequality in \cite{DS} and plays a crucial role in our proof later.

\begin{Proposition}\label{prop:local}
Let $0<4r<\rho\leq 1$. It holds that
\ben\label{eq:revised local energy}
&&k^{-2}A(r)+E(r)\nonumber\\
&&\leq C_0k^4\big(\frac{r}{\rho}\big)^2A(\rho)+C_0k^{-1}\big(\frac{\rho}{r}\big)^{3}\big[C(\rho)+C^{\f 13}(\rho)D^{\f 23}(\rho)]+C_0\big(\frac{\rho}{r}\big)^{2}C^{\f 13}(\rho)F^{\f13}(\rho),
\een
where $1\leq k\leq \frac{\rho}{r}$ and the constant $C_0$ is independent of $k, r, \rho$.
\end{Proposition}

\no{\bf Proof.}\,Let $\zeta$ be a cutoff function, which vanishes
outside of $B_{\rho}$ and equals 1 in $B_{\rho/2}$, satisfying
$$|\nabla\zeta|\leq C_0\rho^{-1},\quad|{\nabla}^2\zeta|\leq C_0\rho^{-2}.$$
Introduce a smooth function as
$$\Gamma(x)=\frac{1}{(k^2r^2+|x|^2)^{2}}, \quad 1\leq k\leq \frac{\rho}{r},$$
which clearly satisfies
$$\triangle\Gamma=\frac{-24k^2r^2}{(k^2r^2+|x|^2)^4}<0.$$
Taking the test function $\phi=\Gamma\zeta$ in the local energy inequality (\ref{eq:local energy}), we obtain that
\beno
&&-\int_{B_{r}^+}|u|^2\zeta\triangle\Gamma dx+2\int_{B_{r}^+}|\nabla u|^2\zeta\Gamma dx\\
&&\leq \int_{B_{\rho}^+}\big[|u|^2(\Gamma\triangle\zeta+2\nabla\Gamma\cdot\nabla\zeta)+u\cdot\nabla\phi(|u|^2+2\pi-2\pi_{B_{\rho}^+})\big]dx +2\int_{B_\rho^+}f u \Gamma\zeta dx
\eeno
It follows from some straightforward computations that
\beno
&&\zeta\Gamma(x,t)\geq C_0(kr)^{-4},\quad-\zeta\triangle\Gamma(x,t)\geq C_0(kr)^{-6}\quad{\rm in} \,\, B_r^+,\\
&&|\nabla\phi|\leq |\nabla\Gamma|\zeta+\Gamma|\nabla\zeta|\leq C_0(kr)^{-5}\quad{\rm in} \,\, B_{\rho}^+,\\
&&|\Gamma\triangle\zeta|+2|\nabla\Gamma\cdot\nabla\zeta|\leq C_0\rho^{-6}\quad{\rm in} \,\, B_{\rho}^+,
\eeno
from which and the H\"{o}lder inequality, one derives that
\beno
 &&k^{-2}A(r)+E(r)\\
 &&\leq C_0k^4\big(\frac{r}{\rho}\big)^2A(\rho)+C_0k^{-1}\big(\frac{\rho}{r}\big)^{3}\rho^{-3}\int_{B_{\rho}^+}(|u|^3+|u||\pi-\pi_{B_{\rho}^+}|)dx\\
 &&+C_0r^{-2}{\rho}^{2}\big(\int_{B_{\rho}^+}|f|^3dx\big)^{\f13}\big(\int_{B_{\rho}^+}|u|^3dx\big)^{\f13}\\
 &&\leq C_0k^4\big(\frac{r}{\rho}\big)^2A(\rho)+C_0k^{-1}\big(\frac{\rho}{r}\big)^{3}\big[C(\rho)+C^{\f 13}(\rho)D^{\f 23}(\rho)]+C_0\big(\frac{\rho}{r}\big)^{2}C^{\f 13}(\rho)F^{\f13}(\rho).
\eeno
The proof is completed. \endproof

\begin{Remark}
If $k=1$ in (\ref{eq:revised local energy}), then we obtain the local energy inequality as the interior case in \cite{DS}:
\ben\label{eq:local-ds}
&&A(r)+E(r)\nonumber\\
&&\leq C_0\big(\frac{r}{\rho}\big)^2A(\rho)+C_0\big(\frac{\rho}{r}\big)^{3}\big[C(\rho)+C^{\f 13}(\rho)D^{\f 23}(\rho)]+C_0\big(\frac{\rho}{r}\big)^{2}C^{\f 13}(\rho)F^{\f13}(\rho).
\een
However, it seems difficult for us to construct an effective iteration scheme based on (\ref{eq:local-ds}) due to the fast decay of the second term in the right hand side of (\ref{eq:local-ds}) (for more details see Remark 4.2). Yet, by choosing a suitable $k$ in Proposition \ref{prop:local}, we can overcome the difficulty, see the proof of Lemma \ref{lem:aed}.
\end{Remark}

\setcounter{equation}{0}
\section{Proof of Theorem \ref{thm:c}}
 In this section, we are going to prove Theorem \ref{thm:c} assuming that Proposition \ref{prop:cd} holds ture, and the next section is devoted to the proof of Proposition \ref{prop:cd}. In fact, due to the Sobolev inequality in Lemma \ref{lem:c1},  it suffices to assume that
\ben\label{eq:C r small}
 C(r)<\varepsilon, \quad {\rm for\,\, any }\quad {0<r<r_0<1}.
\een
First, we need some technical lemmas on the scaling quantities introduced before.
\begin{Lemma}\label{lem:c1}
For any $0<r<r_0$, there holds
\ben\label{eq:C r D r}
C(r)\leq C_0 E(r)^{\f 32},\quad D(r)\leq C_0 D_1^{\f 32}(r).
\een
\end{Lemma}

\no{\bf Proof.}\,These follow from the Sobolev imbedding inequality in the six dimensional space.
\endproof

\begin{Lemma}\label{lem:p1}
For any $0<4r<\rho<r_0$, there holds
\ben\label{eq:D 1}
D_1(r)\leq C_0\big(\f r{\rho}\big)^{3-\frac6p}\big(E^{\f12}(\rho)+D_1(\rho)\big)+C_0\big(\f{\rho}r\big)^2\big(E^{\f12}(\rho)C^{\f13}(\rho)+F^{\f13}(\rho)\big),
\een
where $p\geq 12,$ and $C_0$ depends on $p$.
\end{Lemma}

\no{\bf Proof.}\,
Choose a domain $\widetilde{B}^+$ with a smooth boundary such that $B_{\rho/2}^+\subset\widetilde{B}^+\subset B_{\rho}^+.$ Let $v$ and $\pi_1$ be the unique solution to the following initial boundary value problem for the Stokes system:
\begin{eqnarray}\label{system}
\left\{
  \begin{array}{lll}
    -\Delta v+\nabla{\pi_1}&=&f-u\cdot\nabla u,\quad \hbox{div}v=0\quad {\rm{in}}\quad\widetilde{B}^+, \\
   ~~~~~~~~(\pi_1)_{\widetilde{B}^+}&=&{\f 1{|\widetilde{B}^+|}}\int_{\widetilde{B}^+}\pi_1 dx=0, \\
   ~~~~~~~~~~~~~~v&=&0\quad {\rm on} \quad \partial{\widetilde{B}^+}.
  \end{array}
\right.
\end{eqnarray}
Then by the estimates for the steady Stokes system in Lemma \ref{lem:stokes1}, one can obtain
\ben\label{eq:pi 1}
&&{\f 1{\rho^2}}\parallel v\parallel_{L^{\f65}(\widetilde{B}^+)}+{\f 1{\rho}}\parallel \nabla v\parallel_{L^{\f65}(\widetilde{B}^+)}+\parallel \nabla^2 v\parallel_{L^{\f65}(\widetilde{B}^+)}+{\f 1{\rho}}\parallel \pi_1\parallel_{L^{\f65}(\widetilde{B}^+)}+\parallel \nabla \pi_1\parallel_{L^{\f65}(\widetilde{B}^+)}\nonumber\\
&\leq& C_0\big(\parallel u\cdot\nabla u\parallel_{L^{\f65}(\widetilde{B}^+)}+\parallel f\parallel_{L^{\f65}(\widetilde{B}^+)}\big)\nonumber\\
&\leq& C_0{\rho}^2\big(E^{\f12}(\rho)C^{\f13}(\rho)+F^{\f13}(\rho)\big)
\een
with $C_0$ depending on $B_1^+$.

On the other hand, let
\ben\label{eq:w}
w=u-v,\quad \pi_2=\pi-\pi_{B_{\rho/2}^+}-\pi_1
\een
Then $w,\pi_2$ solve the following boundary value problem:
$$-\Delta w+\nabla{\pi_2}=0,\quad \hbox{div}w=0\quad {\rm in} \quad\widetilde{B}^+,$$
$$w=0\quad {\rm on} \quad \partial{\widetilde{B}^+}\cap\{x_6=0\}.$$
Then the local estimate near the boundary for the steady Stokes system in Lemma \ref{lem:stokes2} for $p\geq 12$ yields that
\beno
\rho^{3-\frac6p}\parallel \nabla \pi_2\parallel_{L^{p}(B_{\rho/4}^+)}
\leq C_0\big({\f 1{\rho^4}}\parallel w\parallel_{L^{\f65}(B_{\rho/2}^+)}+{\f 1{\rho^3}}\parallel \nabla w\parallel_{L^{\f65}(B_{\rho/2}^+)}+{\f 1{\rho^3}}\parallel \pi_2\parallel_{L^{\f65}(B_{\rho/2}^+)}\big).
\eeno
It follows from the Sobolev inequality, (\ref{eq:w}), and (\ref{eq:pi 1}) that
\ben\label{eq:pi 2}
&&\rho^{3-\frac6p}\parallel \nabla \pi_2\parallel_{L^{p}(B_{\rho/4}^+)}\nonumber\\
&\leq&C_0\big({\f 1{\rho^3}}\parallel \nabla w\parallel_{L^{\f65}(B_{\rho/2}^+)}+{\f 1{\rho^3}}\parallel \pi_2\parallel_{L^{\f65}(B_{\rho/2}^+)}\big)\nonumber\\
&\leq&C_0\big({\f 1{\rho^3}}\parallel \nabla u\parallel_{L^{\f65}(B_{\rho/2}^+)}+{\f 1{\rho^3}}\parallel \nabla v\parallel_{L^{\f65}(B_{\rho/2}^+)}+{\f 1{\rho^2}}\parallel \nabla \pi \parallel_{L^{\f65}(B_{\rho/2}^+)}+{\f 1{\rho^3}}\parallel \pi_1\parallel_{L^{\f65}(B_{\rho/2}^+)}\big)\nonumber\\
&\leq& C_0\big(E^{\f12}(\rho)+E^{\f12}(\rho)C^{\f13}(\rho)
+F^{\f13}(\rho)+D_1(\rho)\big).
\een
Combining the two estimates (\ref{eq:pi 1}) and (\ref{eq:pi 2}) shows that
\beno
D_1(r)&=&{\f 1{r^2}}\big(\int_{B_r^+}|\nabla p|^{\f65}dx\big)^{\f56}\\
&\leq&C_0{\f 1{r^2}}\big(\parallel\nabla \pi_1\parallel_{L^{\f65}(B_{r}^+)}+r^{5-\frac6p}\parallel\nabla \pi_2\parallel_{L^{p}(B_{r}^+)}\big)\\
&\leq&C_0\big(\f{\rho}r\big)^2\big(E^{\f12}(\rho)C^{\f13}(\rho)+F^{\f13}(\rho)\big)+C_0\big(\f r{\rho}\big)^{3-\frac6p}\big(E^{\f12}(\rho)+E^{\f12}(\rho)C^{\f13}(\rho)
+F^{\f13}(\rho)+D_1(\rho)\big)\\
&\leq&C_0\big(\f r{\rho}\big)^{3-\frac6p}\big(E^{\f12}(\rho)+D_1(\rho)\big)+C_0\big(\f{\rho}r\big)^2\big(E^{\f12}(\rho)C^{\f13}(\rho)+F^{\f13}(\rho)\big).
\eeno
Thus the proof is completed.
\endproof

{\bf Proof of Theorem \ref{thm:c}:} Take $0<16r<4\rho<\kappa<r_0$ and set
\beno
G(r)=A(r)+E(r)+\varepsilon^{\f16}D_1(r)+\varepsilon^{\f1{12}}F^{\f13}(r).
\eeno
It follows from (\ref{eq:C r D r}) and (\ref{eq:D 1}) for $p=12$ that
\ben\label{eq:D 1 k}
D_1(\rho)
&\leq& C_0\big(\f {\rho}{\kappa}\big)^{\f52}\big(E^{\f12}(\kappa)+D_1(\kappa)\big)+C_0\big(\f{\kappa}{\rho}\big)^2\big(E^{\f12}(\kappa)C^{\f13}(\kappa)+F^{\f13}(\kappa)\big)\nonumber\\
&\leq& C_0\big(\f{\kappa}{\rho}\big)^2\big(E(\kappa)+F^{\f13}(\kappa)\big)+C_0\big(\f {\rho}{\kappa}\big)^{\f52}\big(D_1(\kappa)+1\big).
\een
Then using the local energy inequality in Proposition \ref{prop:local} for $k=1$, (\ref{eq:C r small}) and (\ref{eq:C r D r}), one can obtain that
\beno
G(r)
&\leq&C_0\big(\frac{r}{\rho}\big)^2A(\rho)+C_0\big(\frac{\rho}{r}\big)^{3}\big[C(\rho)+C^{\f 13}(\rho)D^{\f 23}(\rho)]+C_0\big(\frac{\rho}{r}\big)^{2}C^{\f 13}(\rho)F^{\f13}(\rho)\\
&&+C_0\varepsilon^{\f16}D_1(r)+C_0\varepsilon^{\f1{12}}F^{\f13}(r)\\
&\leq&C_0\big(\frac{r}{\rho}\big)^2A(\rho)+C_0\big(\frac{\rho}{r}\big)^{3}\varepsilon^{\f13}D_1(\rho)
+C_0\big(\frac{\rho}{r}\big)^{2}\varepsilon^{\f13}F^{\f13}(\rho)
+C_0\varepsilon^{\f1{12}}F^{\f13}(r)\\
&&+C_0\varepsilon^{\f16}D_1(r)+C_0\big(\frac{\rho}{r}\big)^{3}\varepsilon.
\eeno
Applying the inequality (\ref{eq:D 1 k}) twice yields
\beno
G(r)
&\leq&C_0\big(\frac{\kappa}{\rho}\big)^{4}\big(\frac{r}{\rho}\big)^2G(\kappa)
+C_0\varepsilon^{\f13}\big(\frac{\rho}{r}\big)^{3}\big(\f{\kappa}{\rho}\big)^2\big(E(\kappa)+F^{\f13}(\kappa)\big)\\
&&+C_0\varepsilon^{\f13}\big(\frac{\rho}{r}\big)^{3}\big(\f {\rho}{\kappa}\big)^{\f52}\big(D_1(\kappa)+1\big)
+C_0\varepsilon^{\f13}\big(\frac{\rho}{r}\big)^{2}\big(\frac{\rho}{\kappa}\big)F^{\f13}(\kappa)
+C_0\varepsilon^{\f1{12}}\big(\frac{r}{\kappa}\big)F^{\f13}(\kappa)\\
&&+C_0\varepsilon^{\f16}\big(\f{\kappa}{r}\big)^2\big(E(\kappa)+F^{\f13}(\kappa)\big)
+C_0\varepsilon^{\f16}\big(\f {r}{\kappa}\big)^{\f52}\big(D_1(\kappa)+1\big)+C_0\big(\frac{\rho}{r}\big)^{3}\varepsilon\\
&\leq&C_0\big[\big(\frac{\kappa}{\rho}\big)^{4}\big(\frac{r}{\rho}\big)^2+\big(\frac{r}{\kappa}\big)^{\f52}+\big(\frac{r}{\kappa}\big)]G(\kappa)
+C_0\varepsilon^{\f13}\big(\frac{\kappa}{\rho}\big)^{2}\big(\frac{\rho}{r}\big)^{3}G(\kappa)\\
&&+C_0\varepsilon^{\f14}\big[\big(\frac{\rho}{r}\big)^{3}\big(\frac{\kappa}{\rho}\big)^2+\big(\frac{\rho}{r}\big)^{2}\big(\frac{\rho}{\kappa}\big)\big]G(\kappa)
+C_0\varepsilon^{\f16}\big[\big(\frac{\rho}{r}\big)^{3}\big(\frac{\rho}{\kappa}\big)^{\f52}+\big(\frac{\kappa}{r}\big)^{2}\big]G(\kappa)
\\
&&+C_0\varepsilon^{\f1{12}}\big(\frac{\kappa}{r}\big)^{2}G(\kappa)
+C_0\varepsilon^{\f13}\big(\frac{\rho}{r}\big)^{3}\big(\frac{\rho}{\kappa}\big)^{\f52}
+C_0\varepsilon^{\f16}\big(\frac{r}{\kappa}\big)^{\f52}+C_0\big(\frac{\rho}{r}\big)^{3}\varepsilon.
\eeno
Set $r=\theta^3\rho,\,\rho=\theta\kappa$ with $0<\theta<\f18$. The above inequality yields that
$$
G(r)\leq C_0\big[\theta^{2}+\varepsilon^{\f13}{\theta}^{-11}+\varepsilon^{\f14}{\theta}^{-11}+\varepsilon^{\f16}{\theta}^{-8}+\varepsilon^{\f1{12}}{\theta}^{-8}\big]G(\kappa)+C_0\varepsilon^{\f13}\theta^{-\f{13}{2}}+C_0\varepsilon^{\f16}\theta^{10}+C_0\varepsilon\theta^{-9}.
$$
Thus choosing $\theta$ small at first, and then $\varepsilon$ small,  one can get
\beno
&&C_0\big[\theta^{2}+\varepsilon^{\f13}{\theta}^{-11}+\varepsilon^{\f14}{\theta}^{-11}+\varepsilon^{\f16}{\theta}^{-8}+\varepsilon^{\f1{12}}{\theta}^{-8}\big]\leq \f12,\\
&&\quad C_0\varepsilon^{\f13}\theta^{-\f{13}{2}}+C_0\varepsilon^{\f16}\theta^{10}+C_0\varepsilon\theta^{-9} \leq \varepsilon_1^3,
\eeno
where $\varepsilon_1$ is the constant as in Proposition \ref{prop:cd}. Thus we obtain the following iterative inequality
\ben\label{eq:g k}
G(\theta^{4}\kappa)\leq \frac12 G(\kappa)+\varepsilon_1^3.
\een

On the other hand, it is easy to see that
\beno
G(r_0)\le C_0,
\eeno
where $C_0$ depends on $r_0$, $\|u\|_{H^1(B_1^+)},\|\nabla p\|_{L^{\f65}(B_1^+)}$ and $\|f\|_{L^6(B_1^+)}$.
Then the standard iteration arguments for (\ref{eq:g k}) ensure that there exists $r_1>0$ such that
\ben\label{eq:g k 2}
G(r)\leq 3 \varepsilon_1^3 \quad {\rm for\,\, all}\quad  r\in(0, r_1).
\een

Next, due to Lemma \ref{lem:p1}, (\ref{eq:g k 2}) and $f\in L^6(B^+_1)$, there exists $r_2>0$ such that
\beno
&D_1(r)\leq C_0(\frac{r}{\rho})^2(\varepsilon_1^{3/2}+D_1(\rho))+C_0(\frac{\rho}{r})^2\varepsilon_1^3,\quad 0<r<\rho<r_2<r_1,\\
&F(\rho)\leq \varepsilon_1^3,\quad 0<\rho<r_2<r_1.
\eeno
Then a standard iteration argumemnt yields that there exists $r_3>0$ such that
\beno
D_1(r)+F(r)\leq C_0\varepsilon_1^{3/2} \quad {\rm for}\quad  0<r<r_3<r_2.
\eeno
Combining (\ref{eq:g k 2}) and the above inequality, we get
\beno
C(r)+D_1(r)+F(r)\leq C_0\varepsilon_1^{3/2}\quad {\rm for}\quad  0<r<r_3.
\eeno
We can assume that $C_0\varepsilon_1^{1/2}\leq 1$ without loss of generality. Hence the proof of Theorem \ref{thm:c}
 is completed due to Proposition \ref{prop:cd}.
\endproof

\setcounter{equation}{0}
\section{Proof of Proposition \ref{prop:cd}}

This section is devoted to the proof of Proposition \ref{prop:cd}. To overcome the lower order term $E^{1/2}(\rho)$ from the pressure decomposition, we'll make full use of the revised local energy inequality in Proposition \ref{prop:local} and the decay estimates (\ref{eq:D 1}) in order to derive an effective iteration (for more details, see the following lemma and Remark \ref{rmk:idea}). Finally, we improve the decay of $\nabla u$ by a similar bootstrap argument as in \cite{DS}, and carry out the boundary estimates as in \cite{Se0}.

\begin{Lemma}\label{lem:aed}
Let $\rho>0$ be a positive constant. Then there exists a $\theta_0$, which is suitably small and independent of $\rho,$ such that
\ben\label{eq:iterative}
&&\theta_0^{\f12}A(\theta_0\rho)+E(\theta_0\rho)+\theta_0^{-6+\f{1}{10}}D_1^{2}(\theta_0\rho)\nonumber\\
&\leq&{\f 14}\big[\theta_0^{\f12}A(\rho)+E(\rho)+\theta_0^{-6+\f{1}{10}}D_1^2(\rho)\big]+C_0\big(E^{\f32}(\rho)+E^2(\rho)+F^{\f23}(\rho)\big),
\een
where $C_0$ is a positive constant independent of $\rho$.
\end{Lemma}

\no{\bf Proof.}\, For any $\theta\in(0,\f14],$ Proposition \ref{prop:local} with $r=\theta\rho$ and (\ref{eq:C r D r}) yield that
\ben\label{eq:iterative1}
&&k^{-2}A(\theta\rho)+E(\theta\rho)\nonumber\\
&&\leq C_0k^4\theta^2A(\rho)+C_0k^{-1}\theta^{-3}\big[C(\rho)+C^{\f 13}(\rho)D^{\f 23}(\rho)]+C_0\theta^{-2}C^{\f 13}(\rho)F^{\f13}(\rho),\nonumber\\
&&\leq C_0k^4\theta^2E(\rho)+C_0k^{-1}\theta^{-3}\big[E^{\f32}(\rho)+E^{\f12}(\rho)D_1(\rho)]+C_0\theta^{-2}E^{\f12}(\rho)F^{\f13}(\rho),
 \een

By choosing $p=600$ in Lemma \ref{lem:p1}, one hets
\beno
D_1(\theta\rho)\leq C_0\theta^{3-\frac1{100}}\big(E^{\f12}(\rho)+D_1(\rho)\big)+C_0\theta^{-2}\big(E(\rho)+F^{\f13}(\rho)\big).
\eeno
Hence
\ben\label{eq:iterative2}
D_1^2(\theta\rho)\leq C_0\theta^{6-\frac1{50}}\big(E(\rho)+D_1^2(\rho)\big)+C_0\theta^{-4}\big(E^2(\rho)+F^{\f23}(\rho)\big).
\een

Set
$$G(r)=k^{-2}A(r)+E(r)+\gamma^{-1}D_1^2(r),\quad {\rm and} \,\,r=\theta\rho,$$
with $\gamma>0$ to be decided. It then follows from
(\ref{eq:iterative1}) and (\ref{eq:iterative2}) that
\beno
G(\theta\rho)
&\leq&C_0k^4\theta^2E(\rho)+C_0k^{-1}\theta^{-3}\big[E^{\f32}(\rho)+E^{\f12}(\rho)D_1(\rho)]+C_0\theta^{-2}E^{\f12}(\rho)F^{\f13}(\rho)\\
&&+C_0\gamma^{-1}\theta^{6-\frac1{50}}\big(E(\rho)+D_1^2(\rho)\big)+C_0\gamma^{-1}\theta^{-4}\big(E^2(\rho)+F^{\f23}(\rho)\big)\\
&\leq&C_0k^4\theta^2E(\rho)+C_0k^{-1}\theta^{-3}\big(E^{\f32}(\rho)+\gamma^{\f12}E(\rho)+\gamma^{-\f12}D_1^2(\rho)\big)+C_0\gamma E(\rho)\\
&&+C_0\gamma^{-1}\theta^{6-\frac1{50}}\big(E(\rho)+D_1^2(\rho)\big)+C_0\gamma^{-1}\theta^{-4}\big(E^2(\rho)+F^{\f23}(\rho)\big)\\
&\leq &C_0[k^4\theta^2+k^{-1}\theta^{-3}\gamma^{\f12}+\gamma+\gamma^{-1}\theta^{6-\frac1{50}}+\theta^{6-\frac1{50}}]G(\rho)\\
&&+C_0k^{-1}\theta^{-3}E^{\f32}(\rho)+C_0\gamma^{-1}\theta^{-4}\big(E^2(\rho)+F^{\f23}(\rho)\big)
\eeno
Choosing $k=\theta^{-\f14}$, $\gamma=\theta^{6-\f{1}{10}}$, we get
$$C_0[k^4\theta^2+k^{-1}\theta^{-3}\gamma^{\f12}+\gamma+\gamma^{-1}\theta^{6-\frac1{50}}+\theta^{6-\frac1{50}}]\leq C_0\theta^{\f{4}{50}}.$$
Finally, $\theta$ is chosen so small such that the required inequality (\ref{eq:iterative}) holds.
The lemma is proved.
\endproof
\begin{Remark}\label{rmk:idea}
It should be noted that Proposition \ref{prop:local} is crucial in the proof of Lemma \ref{lem:aed}.
Recall the local energy inequality (\ref{eq:local-ds}) as in \cite{DS}
and the pressure decomposition in Lemma \ref{lem:p1}:
\ben\label{eq:D 1---}
D_1(r)\leq C_0\big(\f r{\rho}\big)^{3-\frac6p}\big(E^{\f12}(\rho)+D_1(\rho)\big)+C_0\big(\f{\rho}r\big)^2\big(E^{\f12}(\rho)C^{\f13}(\rho)+F^{\f13}(\rho)\big),
\een
for $0<4r<\rho<1$ and $p\geq 12$. The two inequalities indicate that $A(r)$, $E(r)$, and $D_1^2(r)$ are the same order, hence it is natural to try to derive an iterative inequality for the following scaled quantity:
$$G(r)=A(r)+E(r)+\gamma^{-1}D_1^2(r),\quad  \,\,r=\theta\rho,$$
with $\gamma>0$ to be decided. It then follows from (\ref{eq:local-ds}) and (\ref{eq:D 1---}) that (as in Lemma 4.1)
\ben\label{eq:G-}
G(\theta\rho)
&\leq &C_0[\theta^2+\theta^{-3}\gamma^{\f12}+\gamma+\gamma^{-1}\theta^{6-\frac{12}{p}}+\theta^{6-\frac{12}{p}}]G(\rho)\nonumber\\
&&+C_0\theta^{-3}E^{\f32}(\rho)+C_0\gamma^{-1}\theta^{-4}\big(E^2(\rho)+F^{\f23}(\rho)\big)
\een
In order to obtain Proposition \ref{prop:cd}, one needs the smallness of
$$C_0[\theta^2+\theta^{-3}\gamma^{\f12}+\gamma+\gamma^{-1}\theta^{6-\frac{12}{p}}+\theta^{6-\frac{12}{p}}],$$
however, which seems to be impossible.

To overcome the difficulty above, we make use of the free parameter in the revised local energy inequality (\ref{eq:revised local energy}) such that the term $\big(\frac{r}{\rho}\big)^2A(\rho)$ in (\ref{eq:local-ds}) plays an important role, and consequently we obtain Lemma 4.1.
\end{Remark}

For $\alpha_0=\log_{\theta_0}^{\f12}$ and  any small $\delta>0,$ we will prove that there exists an integer $m=m(\delta)$ and an increasing sequence of real number $\{\alpha_k\}_{k=1}^m\in(\alpha_0,2)$ with $\alpha_m>2-\delta.$ For these fixed $\delta$ and $m(\delta)$, we obtain the following proposition.
\begin{Proposition}\label{prop:uniform holder}
For $\rho_0\in(0,1)$, there exits a  small constant $\varepsilon_0>0$ satisfying the following property that if
\beno
C(\rho_0)+D_1(\rho_0)+F(\rho_0)\leq\varepsilon_0,
\eeno
then, for any given $\delta>0$, any $0<\rho<\rho_0/8$ and $|x_0|\leq \rho_0/8$ with $x_0=(x',x_6)$ and $x_6=0$, there exists an increasing sequence of real numbers $\{\alpha_k\}_{k=0}^m\in[\alpha_0,2)$ with $\alpha_m>2-\delta$ such that the following inequality holds uniformly
\ben\label{eq:alpha k}
A(\rho,x_0)+E(\rho,x_0)+C^{\f23}(\rho,x_0)+D_1^2(\rho,x_0)\leq C_0{\rho}^{\alpha_k},\quad k=0,1,\cdots,m
\een
where $C_0$ is a positive constant depending on $\delta$ and $k$, but independent of $\rho$.
\end{Proposition}
\no{\bf Proof.}
We will prove this by indiction. Without loss of generality, assume that $x'=0$.
By Proposition \ref{prop:local} and (\ref{eq:C r D r}), it holds that
\beno
&&A(\rho_0/8)+E(\rho_0/8)\\
&\leq&C_0A(\rho_0)+C_0\big[C(\rho_0)+C^{\f 13}(\rho_0)D_1(\rho_0)]+C_0C^{\f 13}(\rho_0) F^{\f13}(\rho_0)\\
&\leq&C_0C^{\f23}(\rho_0)+C_0\big[C(\rho_0)+C^{\f 13}(\rho_0)D_1(\rho_0)]+C_0F^{\f12}(\rho_0)\\
&\leq&C_0(\varepsilon_0+{\varepsilon_0}^{\f12}+\varepsilon_0^{\f23}+\varepsilon_0^{\f43}).
\eeno
Thus we can choose  $\varepsilon_3>0$ first, then $\varepsilon_0=\varepsilon_0(\varepsilon_3)>0$ small enough such that
\ben\label{eq:varepsilon0 small 1}
C_0\big(\varepsilon_3^{\f12}+{\varepsilon_3}\big)<1/4,\quad C_0(\varepsilon_0+{\varepsilon_0}^{\f12}+\varepsilon_0^{\f23}+\varepsilon_0^{\f43})<\varepsilon_3/4,
\een
and
\ben\label{eq:varepsilon0 small 2}\varphi(\rho_0)=\theta_0^{\f{1}{2}}A(\rho_0/8)+E(\rho_0/8)+\theta_0^{-6+\f{1}{10}} D_1^2(\rho_0/8)\leq\varepsilon_3.
\een
It follows from (\ref{eq:varepsilon0 small 1})-(\ref{eq:varepsilon0 small 2}) and (\ref{eq:iterative}) that
\ben\label{eq:phi k}
\varphi(\theta_0^k\rho_0)\leq\varepsilon_3.
\een
Furthermore, by (\ref{eq:phi k}) and (\ref{eq:varepsilon0 small 1}), one has
$$C_0(E(\theta_0^{k-1}\rho_0/8)^{1/2}+E(\theta_0^{k-1}\rho_0/8))<\frac14.$$
Then (\ref{eq:phi k}) and (\ref{eq:iterative}) imply that
\beno
 \varphi(\theta_0^k\rho_0)\leq{\f12}\varphi(\theta_0^{k-1}\rho_0)+C_1(\theta_0^{k-1}\rho_0)^4,
\eeno
where $C_1$ is a constant independent of $k$ and we have used that fact that $f\in L^6$.
One can iterate the above inequality to reach
\beno
\varphi(\theta_0^k\rho_0)&\leq&(\f12)^k\varphi(\rho_0)+C_1\rho_0^4\sum_{j=0}^{k-1}(\f12)^j(\theta_0^{k-1-j})^4\\
&\leq& (\f12)^k\big[\varphi(\rho_0)+\f{2C_1}{1-2\theta_0^4}\rho_0^4\big]
\leq (\f12)^k\big[\varphi(\rho_0)+\f{2C_1}{1-\theta_0}\rho_0^4\big].
\eeno
For $\rho\in(0,\rho_0/8),$ we can find $k$ such that $\theta_0^k{\f{\rho_0}8}<\rho<\theta_0^{k-1}{\f{\rho_0}8}.$ Note that $\alpha_0=\log_{\theta_0}^{\f12},$ then we have
\beno
\theta_0^{\f{1}{2}} A(\rho)+E(\rho)+\theta_0^{-6+\f{1}{10}} D_1^2(\rho)
&\leq& C(\theta_0)\varphi(\theta_0^{k-1}\rho_0)\\
&\leq& C(\theta_0)(\f12)^k\big[\varphi(\rho_0)+\f{2C_1}{1-\theta_0}\rho_0^4\big]\\
&\leq& C(\theta_0)(\theta_0^k)^{\alpha_0}\big[\varphi(\rho_0)+\f{2C_1}{1-\theta_0}\rho_0^4\big]\\
&\leq& C(\theta_0)\rho^{\alpha_0}\big[\varphi(\rho_0)+\f{2C_1}{1-\theta_0}\rho_0^4\big]
\leq C_0\rho^{\alpha_0},
\eeno
where $C_0=C_0(\theta_0,\varphi(\rho_0),C_1,\rho_0)$.  By (\ref{eq:C r D r}) in Lemma \ref{lem:c1}, we can get similar estimates for $C(\rho),$ thus the case $k=0$ is proved.

Assume that (\ref{eq:alpha k}) is true for $m=k$, i.e.
\ben\label{eq:alpha m=k}
A(\rho)+E(\rho)+C^{\f23}(\rho)+D_1^2(\rho)\leq C_0{\rho}^{\alpha_k}.
\een
When $m=k+1,$ we will estimate $A(\rho)+E(\rho).$

Let $\alpha+\beta+\gamma=1,$ with  $\alpha,\beta,\gamma\in(0,1)$ to be decided. By Proposition \ref{prop:local}, (\ref{eq:C r D r}) and the assumption for $m=k$, one can obtain that
\beno
A(\rho)+E(\rho)
&\leq& C_0\rho^{2\alpha}A({\rho}^{\beta+\gamma})+C_0\rho^{-3\alpha}
E({\rho}^{\beta+\gamma})^{\f32}\\
&&+C_0\rho^{-3\alpha}
E({\rho}^{\beta+\gamma})^{\f12}D_1({\rho}^{\beta+\gamma})
+C_0\rho^{-2\alpha}
E({\rho}^{\beta+\gamma})^{\f12}F^{\f13}({\rho}^{\beta+\gamma})\\
&\leq& C_0\rho^{2\alpha}A({\rho}^{\beta+\gamma})+C_0\rho^{-3\alpha}
E({\rho}^{\beta+\gamma})^{\f32}+C_0\rho^{-\f32\alpha}F^{\f12}({\rho}^{\beta+\gamma})\\
&&+C_0\rho^{-3\alpha}
E({\rho}^{\beta+\gamma})^{\f12}D_1({\rho}^{\beta+\gamma}).
\eeno
It follows from (\ref{eq:alpha m=k}) and  $f\in L^6$ that
\ben\label{eq:alpha m=k 1}
&&A(\rho)+E(\rho)\nonumber\\
&\leq& C_0\rho^{2\alpha+(\beta+\gamma)\alpha_k}+C_0\rho^{-3\alpha+\f{3(\beta+\gamma)\alpha_k}2}
+C_0\rho^{-\f32\alpha+3(\beta+\gamma)}
+C_0\rho^{-3\alpha}
E({\rho}^{\beta+\gamma})^{\f12}D_1({\rho}^{\beta+\gamma})\nonumber\\
&\leq& C_0\rho^{\alpha_k+(2-\alpha_k)\alpha}+C_0\rho^{\alpha_k(1+\f\beta2+\f\gamma2-\alpha)-3\alpha}
+C_0\rho^{-\f32\alpha+3(\beta+\gamma)}+C_0\rho^{-3\alpha}E({\rho}^{\beta+\gamma})^{\f12}D_1({\rho}^{\beta+\gamma}).\nonumber\\
\een
For the last term, by Lemma \ref{lem:p1} with $p=12$ and (\ref{eq:alpha m=k}), one can derive that
\ben\label{eq:alpha m=k 2}
\rho^{-3\alpha}E({\rho}^{\beta+\gamma})^{\f12}D_1({\rho}^{\beta+\gamma})
&\leq& C_0\rho^{-3\alpha}{\rho}^{\f{(\beta+\gamma)\alpha_k}2}\big(\rho^{\f{5\beta}{2}}\rho^{\f{\gamma\alpha_k}{2}}+\rho^{-2\beta}\rho^{\gamma\alpha_k}
+\rho^{-2\beta}\rho^{2\gamma}\big)\nonumber\\
&\leq& C_0\rho^{-3\alpha}{\rho}^{\f{(\beta+\gamma)\alpha_k}2}\big(\rho^{\f{5\beta}{2}}\rho^{\f{\gamma\alpha_k}{2}}+\rho^{-2\beta}\rho^{\gamma\alpha_k}\big)\nonumber\\
&\leq& C_0\rho^{\alpha_k(1-\alpha-\f\beta2)+\f52\beta-3\alpha}+C_0\rho^{\alpha_k(1+\f\gamma2-\alpha-\f\beta2)-3\alpha-2\beta}.
\een
Since $\alpha_k<2,$ one can choose
$$\alpha=\f{\alpha_k}{100+5\alpha_k},\quad \beta=\f{4\alpha_k}{100+5\alpha_k},\quad \gamma=\f{100}{100+5\alpha_k}.$$
Now we define  $\alpha_{k+1}$ as :
\beno
\alpha_{k+1}
&=&\min\{{\alpha_k+(2-\alpha_k)\alpha}, {\alpha_k(1+\f\beta2+\f\gamma2-\alpha)-3\alpha}, {-\f32\alpha+3(\beta+\gamma)},\\
&&{\alpha_k(1-\alpha-\f\beta2)+\f52\beta-3\alpha}, {\alpha_k(1+\f\gamma2-\alpha-\f\beta2)-3\alpha-2\beta}\}\\
&=&\min\{\f{102+4\alpha_k}{100+5\alpha_k}\alpha_k,\f{147+6\alpha_k}{100+5\alpha_k}\alpha_k,\f{300+\f{21}2\alpha_k}{100+5\alpha_k},\f{107+2\alpha_k}{100+5\alpha_k}\alpha_k,\f{139+2\alpha_k}{100+5\alpha_k}\alpha_k\}\\
&=&{\f{102+4\alpha_k}{100+5\alpha_k}}\alpha_k\in(\alpha_k,2).
\eeno
Thus it follows from (\ref{eq:alpha m=k 1}) and (\ref{eq:alpha m=k 2}) that
\ben\label{eq:alpha m=k 3}
A(\rho)+E(\rho)+C^{\f23}(\rho)\leq C_0{\rho}^{\alpha_{k+1}},
\een
where $C(\rho)$ is estimated directly by using Lemma \ref{lem:c1}.

Next, we estimate $D_1(\rho).$ Since $\alpha_{k+1}<2,$ by Lemma \ref{lem:p1} with $p=12$ , the inequality
\beno
D_1(\gamma\rho)&\leq& C_0\gamma^{\f52}\big(E^{\f 12}(\rho)+D_1(\rho)\big)+C_0\gamma^{-2}\big(E^{\f 12}(\rho)C^{\f 13}(\rho)+F^{\f 13}(\rho)\big)\\
&\leq&C_0\big(\gamma^{\f52} D_1(\rho)+\gamma^{\f52}\rho^{\f{\alpha_{k+1}}2}+\gamma^{-2}\rho^{\alpha_{k+1}}+\gamma^{-2}\rho^2\big)\\
&\leq&C_0\big(\gamma^{\f52} D_1(\rho)+\gamma^{-2}\rho^{\f{\alpha_{k+1}}2}\big)
\eeno
holds for any $\gamma,\rho\in(0,1).$

For fixed $\rho_1\in(0,\f{\rho_0}{8})$, let $\theta=\f52-\f{\alpha_{k+1}}2$ . Then we can iterate the above inequality to reach
\beno
D_1(\gamma^m\rho_1)
&\leq&C_0\big[\gamma^{\f52} D_1(\gamma^{m-1}\rho_1)+\gamma^{-2}(\gamma^{m-1}\rho_1)^{\f{\alpha_{k+1}}2}\big]\\
&\leq&C_0\big[\gamma^{\f52m}D_1(\rho_1)+\sum_{j=0}^{m-1}\gamma^{-2}\gamma^{\theta j}(\gamma^{m-1}\rho_1)^{\f{\alpha_{k+1}}2}\big]\\
&\leq&C_0\big[(\gamma^{m}\rho_1)^{\f{\alpha_{k+1}}2}\gamma^{(\f52-\f{\alpha_{k+1}}2)m}\rho_1^{\f{\alpha_k-\alpha_{k+1}}2}+\gamma^{-2}\gamma^{-\f{\alpha_{k+1}}2}(\gamma^{m}\rho_1)^{\f{\alpha_{k+1}}2}\sum_{j=0}^{m-1}\gamma^{\theta j}\big]\\
&\leq&C_0\big[(\gamma^m\rho_1)^{\f{\alpha_{k+1}}2}+\gamma^{-3}\f1{1-\gamma^{\theta}}(\gamma^{m}\rho_1)^{\f{\alpha_{k+1}}2}\big]\\
&\leq&C_0(\gamma^m\rho_1)^{\f{\alpha_{k+1}}2},
\eeno
where we have used $\alpha_{k+1}<2$, and $C_0$ is independent of $m$.
This yields  the desired estimate for $D(\rho)$. Combining this with  (\ref{eq:alpha m=k 3}) shows  (\ref{eq:alpha m=k}) for the case $k+1$.  Hence (\ref{eq:alpha k}) is proved.

Finally, it follows from the choice of $\alpha_k,$  that $\alpha_k$ is increasing as $k\rightarrow\infty,$ and $\alpha_k\leq 2$. Thus there exists a limit $\bar{\alpha}$ of  $\alpha_k$ satisfying $0<\bar{\alpha} \leq 2.$ By the definition of $\alpha_{k+1},$ we claim that $\bar{\alpha}=2.$ The proof is completed.
\endproof

{\bf Proof of Proposition \ref{prop:cd}:} It follows from Proposition \ref{prop:uniform holder}  that, for any small $\delta>0$ and $0<\rho<\rho_0/8$, there hold
\ben\label{eq:u holder}
\int_{B_\rho^+}|u(x)|^2&\leq&C_0\rho^{6-\delta},
\een
\ben\label{eq:nabla u holder}
\int_{B_\rho^+}|\nabla u(x)|^2&\leq&C_0\rho^{4-\delta},
\een
\ben\label{eq:u3 holder}
\int_{B_\rho^+}|u(x)|^3dx&\leq&C_0\rho^{6-\f32\delta}.
\een
By (\ref{eq:u holder}), there exists $\rho_1\in(\rho/2,\rho)$ such that
\ben\label{eq:u holder2}
\int_{S_{\rho_1}^+}|u(x)|^2dx\leq C_0{\rho}^{5-\delta},
\een
where $S_{\rho_1}^+=\{x; |x|=\rho_1, x_6\geq 0\}$.

Let $v$ be the unique $H^1$ solution to the Laplace equation
\begin{equation}\nonumber
\left\{\begin{array}{ll}
\displaystyle
\triangle{v}=0,&\hbox{in}\quad{B_{\rho_1}^+},\\
v=0, &\hbox{on}\quad{\partial B_{\rho_1}^+}\cap\{x; x_6=0\},\\
v=u, &\hbox{on}\quad{S_{\rho_1}^+}.
\end{array}\right.
\end{equation}

Now we extend $v(x)$ from $B_{\rho_1}^+$ to $B_{\rho_1}.$ For convenience, we still write $x=(x',x_6)$ and define $\widetilde{v}(x)$ to be
the  odd extension of $v(x)$ from $B_{\rho_1}^+$ to $B_{\rho_1}$  as:
\[\widetilde{v}(x',x_6)=\begin{cases}
v(x',x_6) & x_6\geq0,\\
-v(x',-x_6) & x_6<0.
\end{cases} \]
Moreover, we set
\[\widetilde{u}(x',x_6)=\begin{cases}
u(x',x_6) & x_6\geq0,\\
-u(x',-x_6) & x_6<0.
\end{cases} \]
Obviously, $\widetilde{v}(x)$ satisfies
\begin{equation}\nonumber
\left\{\begin{array}{ll}
\displaystyle
\triangle{\widetilde{v}(x)}=0, & \hbox{in}\quad{B_{\rho_1}},\\
\widetilde{v}=\widetilde{u}, &\hbox{on}\quad{S_{\rho_1}}=\{x; |x|=\rho_1\}.
\end{array}\right.
\end{equation}

Then by the standard estimates for harmonic functions and (\ref{eq:u holder2}), we get
\ben\label{eq:nabla v}
\sup_{B_{{\rho_1}/2}^+}|\nabla v|&\leq& \sup_{B_{\rho_1/2}}|\nabla \widetilde{v}|\leq C_0{\rho_1^{-6}}\int_{S_{\rho_1}}|\widetilde{u}(x)|dx\nonumber\\
&\leq& C_0{\rho_1^{-6}}\int_{S_{\rho_1}^+}|u(x)|dx\leq C_0{\rho}^{-1-\delta/2}.
\een

On the other hand, we let $w=u-v\in H^1(B_{\rho_1}^+),$ then $w$ satisfies the stationary Stokes equation
\begin{equation}\nonumber
\left\{\begin{array}{ll}
\displaystyle
\triangle{w}-\nabla{\pi}=u\cdot\nabla u-f, & \hbox{in}\quad{B_{\rho_1}^+},\\
w=0, &\hbox{on}\quad\partial{B_{\rho_1}^+}.
\end{array}\right.
\end{equation}
Then by the classical $L^p$ estimates for Stokes equations in Lemma \ref{lem:stokes1}, we have
\beno
\|\nabla w\|_{L^{6/5}(B_{\rho_1}^+)}\leq C_0 \rho_1\|u\|_{L^{3}(B_{\rho_1}^+)} \|\nabla u\|_{L^{2}(B_{\rho_1}^+)}+C_0\rho_1\| f\|_{L^{6/5}(B_{\rho_1}^+)}.
\eeno
This, together with the assumption  that $f\in L^6(B_1^+)$ and (\ref{eq:nabla u holder})-(\ref{eq:u3 holder}), yields that
\ben\label{eq:nabla w}
\|\nabla w\|_{L^{6/5}(B_{\rho_1}^+)}
\leq C_0{\rho_1^{5-\delta}}+C_0{\rho_1^{5}}
\leq C_0{\rho_1^{5-\delta}}.
\een
Since $|\nabla u|\leq|\nabla w|+|\nabla v|,$ combining (\ref{eq:nabla v}) and (\ref{eq:nabla w}), for any $r\in(0,\rho/4),$ we obtain that
\beno
\int_{B_r^+}|\nabla u|^{\f{6}{5}}dx &\leq& C_0{\rho^{6-\f65\delta}}+C_0{r^{6}}{\rho^{-\f{6}{5}-\f{3}{5}\delta}}.
\eeno
Taking $r={\f14}{\rho}^{\f{6}{5}-\f1{10}\delta}$, we derive that
\ben\label{eq:nabla u}
\int_{B_r^+}|\nabla u|^{\f{6}{5}}dx &\leq& C_0{r^{\alpha}},
\een
where
$$\alpha=\f{6-\f65\delta}{\f{6}{5}-\f1{10}\delta}>\f{24}5+\frac{1}{10},$$
for a sufficiently small $\delta>0$.

It should be noted that (\ref{eq:nabla u}) and Proposition \ref{prop:dong} imply that
$u$ is  H\"{o}lder continuity in a neighborhood of $x_0$, where $x_0\in B_{\rho_0/8}^+$. In fact, the following arguments are similar to that in \cite{Se0} and we sketch its proof for completeness.

Let $\Gamma=\{x; |x|\leq \rho_0, x_6=0\}$. For  any $x_0\in B_{\rho_0/8}^+$,  we define $r'=dist\{x_0,\Gamma\}=dist\{x_0,x^*\}$ with $x^*\in \Gamma$.

{\bf Case I: $r'\geq \rho_0/{32}$.} Note that $B_{r'}(x_0)\subset B^+_{\rho_0/2}(x^*)$, then one has
\beno
&&r'^{-3}\int_{B_{r'}(x_0)}|u|^3dx+ r'^{-3}\int_{B_{r'}(x_0)}|\pi-\pi_{B_{r'}(x_0)}|^{3/2}dx+r'^{3}\int_{B_{r'}(x_0)}|f|^3dx\\
&&\leq C_0\big(C(\rho_0)+D_1^{3/2}(\rho_0)+F(\rho_0)\big)\leq C_0(\varepsilon+\varepsilon^{3/2})\leq\varepsilon_2.
\eeno
Hence by (\ref{eq:dong}) in Proposition \ref{prop:dong}, we obtain, for any $0<\rho<r'$
\ben\label{ieq:3.10}
\int_{B_{\rho}(x_0)}|\nabla u|^{\f{6}{5}}dx \leq C_0\rho^{\f{24}5+\f{2}{25}}.
\een

{\bf Case II: $r'\leq \rho_0/{32}$.} First,  for $r'/2<\rho< \rho_0/{32}$, since $B_{\rho}(x_0)\subset B^+_{4\rho}(x^*)$, the boundary estimate (\ref{eq:nabla u}) yields that
\ben\label{ieq:3.11}
\int_{B_{\rho}(x_0)\cap B_1^+}|\nabla u|^{\f{6}{5}}dx \leq C_0 \int_{B_{4\rho}^+(x^*)}|\nabla u|^{\f{6}{5}}dx\leq C_0\rho^{\f{24}5+\frac{1}{10}}.
\een
 On the other hand, if $\rho\leq r'/2$, then $B_{r'}(x_0)\subset B^+_{4r'}(x^*)$. By (\ref{eq:alpha k}) in Proposition \ref{prop:uniform holder}, one has
\beno
r'^{-3}\int_{B_{r'}(x_0)}|u|^3dx+ r'^{-3}\int_{B_{r'}(x_0)}|\pi-\pi_{B_{r'}(x_0)}|^{3/2}dx+r'^{3}\int_{B_{r'}(x_0)}|f|^3dx
\leq C_0{r'}<\varepsilon_2,
\eeno
if $\rho_0$ is sufficiently small.
Again, one can use  (\ref{eq:dong}) in Proposition \ref{prop:dong} to derive that
\ben\label{ieq:3.12}
\int_{B_{\rho}(x_0)}|\nabla u|^{\f{6}{5}}dx \leq C_0{\rho^{\f{24}5+\f{2}{25}}},\quad \forall \,\,0<\rho\leq r'/2.
\een

Combining the inequalities (\ref{eq:nabla u})-(\ref{ieq:3.12}) shows that
 for any $x_0\in B_{\rho_0/8}^+$ and $0<\rho<\rho_0/8$,
\ben\label{ieq:3.13}
\int_{B_{\rho}(x_0)\cap B_1^+}|\nabla u|^{\f{6}{5}}dx \leq C_0\rho^{\f{24}5+\f{2}{25}}.
\een
Then the Morrey lemma yields that $u$ is H\"{o}lder continuity in a neighborhood of $0$. Hence we complete the proof of Proposition \ref{prop:cd}.
\endproof

\setcounter{equation}{0}
\section{Boundary regularity of the 6D steady MHD equations }
In this section, we will extend the previous boundary regularity results to the 6-D steady-state incompressible Magneto-hydrodynamics equations (\ref{eq:SMHD}).
As introduced in Section 2, the scaled invariant quantities $D(r)$, $D_1(r)$ and $F(r)$ are invariable, and  we also need some new ones involving both velocity field $u$ and magnetic field $b$ as follows:
$$A(u,b;r)=r^{-4}\int_{B_r^+}|u(x)|^2+|b(x)|^2dx,\quad C(u,b;r)=r^{-3}\int_{B_r^+}|u(x)|^3+|b(x)|^3,$$
$$E(u,b;r)=r^{-2}\int_{B_r^+}|\nabla u(x)|^2+|\nabla b(x)|^2dx.$$
Moreover, $A(u,0;r)=A(u,r)$, $A(0,b;r)=A(b,r)$, and similar notations hold for $C(u,b;r)$ and $E(u,b;r)$.

The following revised local energy inequality is similar to Proposition \ref{prop:local}.
\begin{Proposition}\label{prop:local-MHD}
Let $0<4r<\rho\leq 1$. It holds that
\ben\label{eq:revised local energy-MHD}
&&k^{-2}A(u,b;r)+E(u,b;r)\nonumber\\
&\leq& C_0k^4\big(\frac{r}{\rho}\big)^2A(u,b;\rho)+C_0k^{-1}\big(\frac{\rho}{r}\big)^{3}\big[C(u,\rho)+C^{\f 13}(u,\rho)D^{\f 23}(\rho)+C^{\f 13}(u,\rho)C^{\f 23}(b,\rho)]\nonumber\\
&&+C_0\big(\frac{\rho}{r}\big)^{2}C^{\f 13}(u,\rho)F^{\f13}(\rho),
\een
where $1\leq k\leq \frac{\rho}{r}$ and the constant $C_0$ is independent of $k, r, \rho$.
\end{Proposition}

Compared with Lemma \ref{lem:p1}, the pressure estimates can be modified as follows:
\begin{Lemma}\label{lem:p1-MHD}
For any $0<4r<\rho<r_0$, there holds
\ben\label{eq:D 1-MHD}
D_1(r)\leq C_0\big(\f r{\rho}\big)^{3-\frac6p}\big(E^{\f12}(u,\rho)+D_1(\rho)\big)+C_0\big(\f{\rho}r\big)^2\big(E^{\f12}(u,b;\rho)C^{\f13}(u,b;\rho)+F^{\f13}(\rho)\big),
\een
where $p\geq 12,$ and $C_0$ depends on $p$.
\end{Lemma}

Considering the velocity field and magnetic field  together, immediately,
we get the $\varepsilon$ regularity criterion involving $u$ and $b$ in the same way as Proposition \ref{prop:cd}.
\begin{Proposition}\label{prop:cd-MHD}
Let $(u,b,\pi,)$ be a suitable weak solution to (\ref{eq:SMHD}) in $B_1^+$ near the boundary $\{x\in B_1, x_6=0\}$.
If there exists ${\rho}_0>0$ and a small positive constant $\varepsilon_1$
such that
\beno
\rho_0^{-3}\|u\|_{L^3(B_{\rho_0}^+)}^3+\rho_0^{-3}\|b\|_{L^3(B_{\rho_0}^+)}^3+\rho_0^{-2}\|\nabla \pi\|_{L^{6/5}(B_{\rho_0}^+)}+\rho_0^3\|f\|_{L^3(B_{\rho_0}^+)}^3<\varepsilon_1
\eeno
Then $0$ is a regular point of $(u,b)$.
\end{Proposition}

{\bf Proof of Theorem \ref{thm:MHD}:} The  proof is similar to that of Theorem \ref{thm:c} with slight
changes, and we'll construct an effective iterative scheme.


let $0<16r<4\rho<\kappa<r_0$. Applying Lemma \ref{lem:c1} and putting $p=12$ in (\ref{eq:D 1-MHD}), one has
\ben\label{eq:5.4}
D_1(\rho)
&\leq& C_0\big(\f {\rho}{\kappa}\big)^{\f52}\big(E^{\f12}(u,\kappa)+D_1(\kappa)\big)+C_0\big(\f{\kappa}{\rho}\big)^2\big(E^{\f12}(u,b;\kappa)C^{\f13}(u,b;\kappa)+F^{\f13}(\kappa)\big)\notag\\
&\leq& C_0\big(\f{\kappa}{\rho}\big)^2\big(E(u,b;\kappa)+A(u,b;\kappa)+F^{\f13}(\kappa)\big)+C_0\big(\f {\rho}{\kappa}\big)^{\f52}\big(D_1(\kappa)+1\big).
\een

Set
\beno
H(r)=A(u,b;r)+E(u,b;r)+\varepsilon^{\f16}D_1(r)+\varepsilon^{\f1{12}}F^{\f13}(r),
\eeno
then it follows from (\ref{eq:revised local energy-MHD}) for $k=1$ and  Lemma \ref{lem:c1}  that
\beno
H(r)
&\leq&C_0\big(\frac{r}{\rho}\big)^2A(u,b;\rho)+C_0\big(\frac{\rho}{r}\big)^{3}\big[C(u,\rho)+C^{\f 13}(u,\rho)(C^{\f 23}(b,\rho)+D^{\f 23}(\rho))]\\
&&+C_0\big(\frac{\rho}{r}\big)^{2}C^{\f 13}(u,\rho)F^{\f13}(\rho)+C_0\varepsilon^{\f16}D_1(r)+C_0\varepsilon^{\f1{12}}F^{\f13}(r)\\
&\leq&C_0\big(\frac{r}{\rho}\big)^2A(u,b;\rho)+C_0\big(\frac{\rho}{r}\big)^{3}\varepsilon^{\f13}E(u,b;\rho)+
C_0\big(\frac{\rho}{r}\big)^{3}\varepsilon^{\f13}D_1(\rho)\\
&&+C_0\big(\frac{\rho}{r}\big)^{2}\varepsilon^{\f13}F^{\f13}(\rho)+C_0\varepsilon^{\f16}D_1(r)+C_0\varepsilon^{\f1{12}}F^{\f13}(r).
\eeno
Moreover, using (\ref{eq:5.4}), we get
\beno
H(r)&\leq&C_0\big[\big(\frac{\kappa}{\rho}\big)^{4}\big(\frac{r}{\rho}\big)^2+\big(\frac{\rho}{r}\big)^{3}\big(\f{\kappa}{\rho}\big)^2\varepsilon^{\f13}\big]H(\kappa)
\\
&&+C_0\varepsilon^{\f13}\big(\frac{\rho}{r}\big)^{3}\big(\f{\kappa}{\rho}\big)^2\big(E(u,b;\kappa)+A(u,b;\kappa)+F^{\f13}(\kappa)\big)+C_0\varepsilon^{\f13}\big(\frac{\rho}{r}\big)^{3}\big(\f {\rho}{\kappa}\big)^{\f52}\big(D_1(\kappa)+1\big)\\
&&+C_0\varepsilon^{\f16}D_1(r)+C_0\varepsilon^{\f13}\big(\frac{\rho}{r}\big)^{2}\big(\frac{\rho}{\kappa}\big)F^{\f13}(\kappa)+C_0\varepsilon^{\f1{12}}\big(\frac{r}{\kappa}\big)F^{\f13}(\kappa)\\
&\leq&C_0\big[\big(\frac{\kappa}{\rho}\big)^{4}\big(\frac{r}{\rho}\big)^2+\big(\frac{\rho}{r}\big)^{3}\big(\f{\kappa}{\rho}\big)^2\varepsilon^{\f13}\big]H(\kappa)
+C_0\varepsilon^{\f13}\big(\frac{\rho}{r}\big)^{3}\big(\f{\kappa}{\rho}\big)^2\big(F^{\f13}(\kappa)\big)\\
&&+C_0\varepsilon^{\f13}\big(\frac{\rho}{r}\big)^{3}\big(\f {\rho}{\kappa}\big)^{\f52}\big(D_1(\kappa)+1\big)
+C_0\varepsilon^{\f1{12}}\big(\frac{r}{\kappa}\big)F^{\f13}(\kappa)\\
&&+C_0\varepsilon^{\f16}\big(\f{\kappa}{r}\big)^2\big(E(u,b;\kappa)+A(u,b;\kappa)+F^{\f13}(\kappa)\big)+C_0\varepsilon^{\f16}\big(\f {r}{\kappa}\big)^{\f52}\big(D_1(\kappa)+1\big)\\
&\leq&C_0\big[\big(\frac{\kappa}{\rho}\big)^{4}\big(\frac{r}{\rho}\big)^2+\big(\frac{r}{\kappa}\big)^{\f52}+\big(\frac{r}{\kappa}\big)]H(\kappa)
+C_0\varepsilon^{\f13}\big(\frac{\kappa}{\rho}\big)^{2}\big(\frac{\rho}{r}\big)^{3}H(\kappa)\\
&&+C_0\varepsilon^{\f14}\big(\frac{\rho}{r}\big)^{3}\big(\frac{\kappa}{\rho}\big)^2H(\kappa)
+C_0\varepsilon^{\f16}\big[\big(\frac{\rho}{r}\big)^{3}\big(\frac{\rho}{\kappa}\big)^{\f52}+\big(\frac{\kappa}{r}\big)^{2}\big]H(\kappa)
\\
&&+C_0\varepsilon^{\f1{12}}\big(\frac{\kappa}{r}\big)^{2}H(\kappa)
+C_0\varepsilon^{\f13}\big(\frac{\rho}{r}\big)^{3}\big(\frac{\rho}{\kappa}\big)^{\f52}
+C_0\varepsilon^{\f16}\big(\frac{r}{\kappa}\big)^{\f52},
\eeno
which follows from the same arguments as (\ref{eq:g k 2}) that
\ben\label{eq:mg k 2}
H(r)\leq 3 \varepsilon_1^3 \quad {\rm for\,\, all}\quad  r\in(0, r_1),
\een
where $\varepsilon_1$ is the constant in Proposition \ref{prop:cd-MHD}.

Now we can employ (\ref{eq:5.4}), $f\in L^6(B^+_1)$ and the standard iterative to derive the existence of $r_2>0$ such that
\beno
D_1(r)+F(r)\leq C_0{\varepsilon_1}^{\f32} \quad {\rm for}\quad  0<r<r_2<r_1,
\eeno
which with (\ref{eq:mg k 2}) yield that
\beno
C(u,b;r)+D_1(r)+F(r)\leq C_0{\varepsilon_1}^{\f32}\quad {\rm for}\quad  0<r<r_2.
\eeno

At last, by choosing $C_0\varepsilon_1^{1/2}\leq 1$ and applying Proposition \ref{prop:cd-MHD}, the proof is completed.
\endproof

\bigskip

\noindent {\bf Acknowledgments.}
The authors would like to thank Professor Zhouping Xin for his valuable discussions and constant encouragements. Wang is supported by NSFC 11301048 and "the Fundamental Research Funds for the Central Universities". Liu is supported in part by the CNPq grant $\sharp$. 501376/2013-1.


\end{document}